\documentclass[12pt]{amsart}
\usepackage{amssymb}
\usepackage{hyperref}
\usepackage{tabularx}
\usepackage{booktabs}
\usepackage{caption}

\newtheorem{thm}{Theorem}[section]
\newtheorem{lem}[thm]{Lemma}
\newtheorem{cor}[thm]{Corollary}
\newtheorem{prop}[thm]{Proposition}
\newtheorem{ex}[thm]{Example}

\newtheorem{con}[thm]{Conjecture}
\newtheorem*{prob*}{Open problem}

\theoremstyle{definition}

\newtheorem{defi}[thm]{Definition}

\theoremstyle{remark}

\newtheorem{rem}[thm]{Remark}
\newtheorem*{rem*}{Remark}

\DeclareMathOperator{\id}{id}

\newcommand{\kringel}{\mathbin{\raise1pt\hbox{$\scriptstyle\circ$}}}
\newcommand{\pkt}{\mathbin{\raise0pt\hbox{$\scriptstyle\bullet$}}}

\newcommand{\C}{\mathbb{C}}

\newcommand{\F}{\mathbb{F}}
\newcommand{\N}{\mathbb{N}}

\newcommand{\ad}{{\rm ad}}

\newcommand{\End}{{\rm End}}
\newcommand{\Der}{{\rm Der}}
\newcommand{\Inn}{{\rm Inn}}
\newcommand{\Out}{{\rm Out}}

\newcommand{\La}{\mathfrak{a}}
\newcommand{\Lb}{\mathfrak{b}}

\newcommand{\Le}{\mathfrak{e}}
\newcommand{\Lf}{\mathfrak{f}}
\newcommand{\Lg}{\mathfrak{g}}
\newcommand{\Lh}{\mathfrak{h}}

\newcommand{\Lk}{\mathfrak{k}}

\newcommand{\Ln}{\mathfrak{n}}

\newcommand{\Ls}{\mathfrak{s}}

\newcommand{\CM}{\mathcal{M}}

\newcommand{\al}{\alpha}
\newcommand{\be}{\beta}

\newcommand{\om}{\omega}

\newcommand{\ov}{\overline}

\newcommand{\ra}{\rightarrow}

\renewcommand{\phi}{\varphi}

\begin{document}


\title[CPA-structures]{A Zassenhaus conjecture and CPA-structures on simple modular Lie algebras}

\author[D. Burde]{Dietrich Burde}
\author[W. Moens]{Wolfgang Alexander Moens}
\address{Fakult\"at f\"ur Mathematik\\
Universit\"at Wien\\
  Oskar-Morgenstern-Platz 1\\
  1090 Wien \\
  Austria}
\email{dietrich.burde@univie.ac.at}
\address{Fakult\"at f\"ur Mathematik\\
Universit\"at Wien\\
  Oskar-Morgenstern-Platz 1\\
  1090 Wien \\
  Austria}
\email{wolfgang.moens@univie.ac.at}

\date{\today}

\subjclass[2000]{Primary 17B30, 17D25}
\keywords{Simple modular Lie algebra, Commutative Post-Lie algebra, CPA structure}

\begin{abstract}
Commutative post-Lie algebra structures on Lie algebras, in short CPA structures, have been studied over fields
of characteristic zero, in particular for real and complex numbers motivated by geometry. A perfect Lie algebra in 
characteristic zero only admits the trivial CPA-structure. In this article we study these structures over fields of 
characteristic $p>0$. We show that every perfect modular Lie algebra in characteristic $p>2$ having a solvable outer derivation 
algebra admits only the trivial CPA-structure. This involves a conjecture by Hans Zassenhaus, saying that the 
outer derivation algebra $\Out(\Lg)$ of a simple modular Lie algebra $\Lg$ is solvable. We try to summarize the known results
on the Zassenhaus conjecture and prove some new results using the classification of simple modular 
Lie algebras by Premet and Strade for algebraically closed fields of characteristic $p>3$. As a corollary we obtain that
that every central simple modular Lie algebra of characteristic $p>3$ admits only the trivial CPA-structure.
 \end{abstract}

\maketitle

\section{Introduction}

Post-Lie algebras and post-Lie algebra structures arise in many areas of mathematics and physics. We have studied these 
structures in the context of geometric structures on Lie groups, in particular for nil-affine actions on Lie groups. 
Here post-Lie algebras arise as a natural common generalization of pre-Lie algebras \cite{HEL,KIM,SEG,BU5,BU19,BU24}
and LR-algebras \cite{BU34, BU38}. We have obtained several results on the existence of post-Lie algebra structures 
\cite{BU33,BU41,BU44,BU51}. Originally post-Lie algebras have been introduced by Vallette \cite{VAL} in connection with 
the homology of partition posets and the study of Koszul operads and by Loday \cite{LOD} within the context of algebraic 
operad triples. Then several authors have studied them in connection with modified Yang-Baxter equations, double Lie algebras, 
$R$-matrices, isospectral flows, Lie-Butcher series and other topics \cite{ELM}. \\
A particular interesting class of post-Lie algebra structures are {\em commutative} structures, so-called {\em CPA-structures}. 
We have studied CPA-structures for several classes of finite-dimensional Lie algebras over a field of characteristic zero,
e.g. for semisimple, perfect and complete Lie algebras, see in \cite{BU51, BU52}. In particular, we have shown that CPA-structures 
on perfect Lie algebras in characteristic zero are trivial. We also extended this result to some classes of perfect
in\-fi\-ni\-te-dimensional Lie algebras \cite{BU60}. The aim of this article is to generalize the result to 
modular Lie algebras over a field of characteristic $p>0$. In particular, we are interested in CPA-structures on simple 
modular Lie algebras. It turns out that we can determine such structures if we know that the outer derivation algebra 
$\Out(\Lg)$ is solvable. But this leads us directly to a conjecture by Hans Zassenhaus, which says that $\Out(\Lg)$ is 
solvable for all simple modular Lie algebras. Zassenhaus posed this conjecture in $1939$ in his book \cite{ZAS}. 
We collect the known results on this conjecture from the literature. For algebraically closed fields of characteristic 
$p>3$ we can use the classification results by Wilson and Strade to verify the conjecture in this case. We obtain that
$\Out(\Lg)$ is solvable of derived length at most three. For $p=2$ and $p=3$ there are counterexamples arising from the 
simple constituents of the classical Lie algebras $\Lf_4$ and $\Lg_2$ respectively. \\
In order to apply the results on the Zassenhaus conjecture to CPA-structures we need another property concerning the 
outer derivation algebra for a Lie algebra $\Lg$. We say that $\Lg$ is {\em taut} if every Lie algebra homomorphism
$\phi\colon \Lg\ra \Der(\Lg)$ maps $\Lg$ into $\Inn(\Lg)$. A simple modular Lie algebra $\Lg$ having a solvable outer
derivation algebra is taut, but the converse need not be true. In fact, the simple modular Lie algebra 
$\mathfrak{psl}_3(\F_3)$ is taut, but its outer derivation algebra is not solvable. This is one of the 
counterexamples to the Zassenhaus conjecture mentioned above. We can try to consider a weaker conjecture here. 
Namely, we conjecture that every simple modular Lie algebra is taut. \\
Using the results on taut Lie algebras and the Zassenhaus conjecture we prove that every CPA-structure on a central
simple Lie algebra of characteristic $p>3$ is trivial. Furthermore we show that every CPA-structure on a taut Lie
algebra in characteristic §$p\neq 2$ is trivial, and that every CPA-structure on a taut Lie
algebra in characteristic §$p= 2$ is either trivial or adjoint.

\section{Preliminaries}

Let $\Lg$ denote a Lie algebra over a field $K$ of arbitrary characteristic. We always assume that $\Lg$ is 
finite-dimensional over $K$. The {\it derived series} of $\Lg$ is defined by
\[
\Lg^{(0)}=\Lg \supseteq \Lg^{(1)} \supseteq \Lg^{(2)} \supseteq \Lg^{(3)}  \supseteq \cdots
\]
where the ideals $\Lg^{(i)}$ are defined by $\Lg^{(i)}=[\Lg^{(i-1)},\Lg^{(i-1)}]$
for all $i\ge 1$. A Lie algebra $\Lg$ is called $k$-step solvable if $\Lg^{(k)}=0$, but
$\Lg^{(k-1)}\neq 0$ for some $k\ge 1$. In this case we call $k$ the {\it solvability class} of $\Lg$, and
denote it by $d(\Lg)$. It is also called the {\it derived length} of $\Lg$. \\
The {\it lower central series} of $\Lg$ is given by
\[
\Lg^0=\Lg \supseteq \Lg^1 \supseteq \Lg^2 \supseteq \Lg^3  \supseteq \cdots
\]
where the ideals $\Lg^i$ are recursively defined by $\Lg^i=[\Lg,\Lg^{i-1}]$
for all $i\ge 1$.

\begin{defi}
Let $\Lg$ be a Lie algebra. The intersection of all ideals $\Lg^{(i)}$ is denoted by
\[
\Lg^{(\infty)}=\bigcap_{i\in \N}\Lg^{(i)}.
\]
Then $\Lg^{(\infty)}=0$ if and only if $\Lg$ is solvable. Similarly we define
\[
\Lg^{\infty}=\bigcap_{i\in \N}\Lg^{i}.
\]
\end{defi}

Let $\Der(\Lg)$ be the Lie algebra of derivations of $\Lg$ and $\Inn(\Lg)$ the ideal 
of inner derivations of $\Lg$. The quotient Lie algebra $\Der(\Lg)/\Inn(\Lg)$ is called {\em outer derivation algebra 
of $\Lg$} and is denoted by $\Out(\Lg)$. In terms of Lie algebra cohomology we have $Z^1(\Lg,\Lg)=\Der(\Lg)$, 
$B^1(\Lg,\Lg)=\Inn(\Lg)$ and $H^1(\Lg,\Lg)=\Out(\Lg)$, where $\Lg$ denotes the adjoint $\Lg$-module. \\[0.2cm]
The {\em centroid} $C(\Lg)$ of a Lie algebra $\Lg$ is the space of $K$-linear maps $\phi\colon \Lg\to \Lg$ commuting 
with all inner derivations of $\Lg$, i.e., satisfying
\[
\phi([x,y])=[\phi(x),y]
\]
for all $x,y\in \Lg$. The centroid of a simple non-associative algebra $A$ is a field, see \cite{JAC}, chapter $10$.
Such an algebra $A$ is then called {\em central simple}, if the centroid $C(A)$ coincides with the base field.
So we have the following definition for Lie algebras.

\begin{defi}
A simple Lie algebra $\Lg$ over $K$ is called {\em central simple}, if $C(\Lg)=K$. 
\end{defi}

A central simple Lie algebra $\Lg$ remains simple for any field extension $F$ over $K$, i.e.,
the Lie algebra $\Lg_F=\Lg\otimes_K F$ is still simple over $F$. In particular, $\Lg$ is still simple over the 
algebraic closure of $K$.

\begin{ex}\label{2.3}
The simple Lie algebra $W(1,\underline{2})^{(1)}$ of dimension $3$  over $\F_2$ is central simple. However,
the simple Lie algebra $W(1,\underline{2})^{(1)}\otimes_{\F_2}\F_4$ of dimension $6$ over $\F_2$ is simple
but not central simple.
\end{ex}

Here $W(1,\underline{2})^{(1)}$ is a Lie algebra of Witt type. The Lie brackets with respect to a $(e,f,h)$
are given by
\[
[e,f]=h,\; [e,h]=e,\; [f,h]=f.
\]
It is up to isomorphism the only simple Lie algebra of dimension $3$ in characteristic $2$ \cite{HIS}. Both
$W(1,\underline{2})^{(1)}$ and $W(1,\underline{2})^{(1)}\otimes_{\F_2}\F_4$ are discussed in \cite{VAU}. \\[0.2cm]
We recall the definition of post-Lie algebra structures on pairs of Lie algebras $(\Lg,\Ln)$ over $K$ \cite{BU41}:

\begin{defi}\label{pls}
Let $\Lg=(V, [\, ,])$ and $\Ln=(V, \{\, ,\})$ be two Lie brackets on a vector space $V$ over 
$K$. A {\it post-Lie algebra structure}, or {\em PA-structure} on the pair $(\Lg,\Ln)$ is a 
$K$-bilinear product $x\cdot y$ satisfying the identities:
\begin{align}
x\cdot y -y\cdot x & = [x,y]-\{x,y\} \label{post1}\\
[x,y]\cdot z & = x\cdot (y\cdot z) -y\cdot (x\cdot z) \label{post2}\\
x\cdot \{y,z\} & = \{x\cdot y,z\}+\{y,x\cdot z\} \label{post3}
\end{align}
for all $x,y,z \in V$.
\end{defi}

Define by  $L(x)(y)=x\cdot y$ and $R(x)(y)=y\cdot x$ the left respectively right multiplication 
operators of the algebra $A=(V,\cdot)$. By \eqref{post3}, all $L(x)$ are derivations of the Lie 
algebra $(V,\{,\})$. Moreover, by \eqref{post2}, the left multiplication
\[
L\colon \Lg\ra \Der(\Ln)\subseteq \End (V),\; x\mapsto L(x)
\]
is a linear representation of $\Lg$. A particular case of a post-Lie algebra structure arises 
if the algebra $A=(V,\cdot)$ is {\it commutative}, i.e., if $x\cdot y=y\cdot x$ is satisfied for all 
$x,y\in V$, so that we have $L(x)=R(x)$ for all $x\in V$. Then the two Lie brackets $[x,y]=\{x,y\}$ 
coincide, and we obtain a commutative algebra structure on $V$ associated with only one Lie 
algebra \cite{BU51}:

\begin{defi}\label{cpa}
A {\it commutative post-Lie algebra structure}, or {\em CPA-structure} on a Lie algebra $\Lg$ 
is a $K$-bilinear product $x\cdot y$ satisfying the identities:
\begin{align}
x\cdot y & =y\cdot x \label{com4}\\
[x,y]\cdot z & = x\cdot (y\cdot z) -y\cdot (x\cdot z)\label{com5} \\
x\cdot [y,z] & = [x\cdot y,z]+[y,x\cdot z] \label{com6}
\end{align}
for all $x,y,z \in V$. 
\end{defi}

\section{A conjecture by Zassenhaus}

Hans Zassenhaus posed in $1939$ in his book \cite{ZAS} the following conjecture.

\begin{con}[Zassenhaus]
The outer derivation algebra of a simple non-abelian Lie algebra in prime characteristic
is solvable.
\end{con}

This is stated between ``Satz $7$'' and ``Satz $8$'' on page $80$ in \cite{ZAS}, where he discusses two other
related conjectures on simple modular Lie algebras. In characteristic zero the outer derivation algebra is trivial 
as a consequence of the first Whitehead Lemma: we have $H^1(\Lg,\Lg)=0$ for a simple Lie algebra $\Lg$ 
of characteristic zero. In prime characteristic Whitehead's Lemma does not hold in general and the outer derivation algebra 
of a simple Lie algebra need not be trivial. 

\begin{rem}
The Zassenhaus conjecture for Lie algebras is an analogue of the {\em Schreier conjecture} for finite groups.
The Schreier conjecture asserts that the outer automorphism group of every finite simple non-abelian group is solvable. 
It was proposed by Otto Schreier in $1926$ and is known to be true as a result of the classification of finite 
simple groups. Up to now no simpler proof is known for it. The classification yields an even stronger statement, namely
that the outer automorphism group of every finite simple non-abelian group is {\em metabelian}.
\end{rem}

What is known about the Zassenhaus conjecture? There are many different results in the literature, in particular
in the context of the classification of simple modular Lie algebras over an algebraically closed field of characteristic
$p>3$. However, there seems to be no general survey available on the outer derivation algebra. Hence we will try to collect 
the relevant results from the literature ourselves. \\
Shortly after Zassenhaus had posed the conjecture, Jacobson published counterexamples in \cite{JA2}, p. $515$,
over an imperfect field of arbitrary characteristic $p>0$. However, it seems that the examples have an obvious
non-trivial ideal, so that the Lie algebras are not simple over the given field.
In general it is reasonable to restrict the Zassenhaus conjecture to {\em algebraically closed fields}. 
Nevertheless we also have results in some cases for finite fields 
or even for arbitrary fields of characteristic $p>0$. In order to use the classification of simple modular Lie algebras 
we need to assume that $p>3$. But again, we also have some results for $p=2$ or $p=3$. \\
According to the classification theorem of Premet and Strade every simple modular Lie
algebra over an algebraically closed field of characteristic $p>3$ belongs to one and only one type: classical, Cartan or 
Melikian.

\subsection{Classical type}
The first class of simple modular Lie algebras consists of the {\em classical Lie algebras} of
Chevalley type. Here it is usually assumed that $p>3$, but the construction works for any field 
of characteristic $p>0$.
Given a simple Lie algebra $L$ over $\C$, Chevalley constructed a basis with integral structure constants. 
Hence we obtain from $L$ a modular Lie algebra $L_F$ by tensoring with any field $F$ of characteristic $p>0$. 
The Lie algebra $L_F$ may fail to be simple, even for $p>3$. But for $p>3$ this happens only in one case. 
In fact, $L_F$ is always simple for $p>3$ except for type $A_n$ with $n\equiv -1 \bmod p$, where $L_F$ has a 
$1$-dimensional center $Z$. Then the quotient $L_F/Z=\mathfrak{psl}_n(F)$ is simple. 
The list of classical simple Lie algebras over $F$ for $p>3$ is given by
\[
A_n, p\nmid n+1,\; \mathfrak{psl}(n+1), p\mid n+1\;B_n,\;C_n,\; D_n,\; G_2,\;F_4,\;E_6,\;E_7,\;E_8. 
\]
Note that the exceptional types here are included in the class of classical simple Lie algebras. \\[0.2cm]
We recall the following result, which is Lemma $2.7$ in \cite{BGP}.

\begin{prop}\label{3.3}
Let $\Lg$ be a classical simple Lie algebra over an algebraically closed field $\F$ of characteristic 
$p>3$. Then $\Out(\Lg)=0$ unless $\Lg=\mathfrak{psl}_{n+1}(\F)$ with $p\mid n+1$ in which case 
$\Der(\Lg)\cong \mathfrak{pgl}_{n+1}(\F)$ and $\Out(\Lg)\cong \F$.
\end{prop}

What can we say in characteristic $p=3$ about the Zassenhaus conjecture for simple modular Lie algebras 
of classical type? Taking the list for $p>3$ it is well known that
all types given there remain simple for $p=3$ except for $\Lg_2$ and $\Le_6$, see \cite{DIE}, p. $204$ or
\cite{STR1}, p. $209$. 
The Lie algebra $\Lg_2$ is no longer simple, but has a simple ideal $I\cong \mathfrak{psl}_3(\F)$, generated
by the short roots, with $\Lg_2/I\cong I$. The Lie algebra $\Le_6$ of dimension $78$ has a $1$-dimensional 
center $Z$, and the quotient $\Le_6/Z$ is a simple Lie algebra of dimension $77$. \\
By Corollary $3.3$ and $3.4$ in \cite{ELK} and by the table in section $4$ of \cite{BBG} we have the following result. 

\begin{prop}
Let $\F$ be an algebraically closed field of characteristic $p=3$ and $\Lg$ be a simple Lie algebra of 
classical type different from $\mathfrak{psl}_{3m}(\F)$, $\Lg_2$, or $\Le_6$. Then $\Out(\Lg)=0$.
\end{prop}

Now the simple Lie algebra $\Lg=\mathfrak{psl}_3(\F)$ immediately yields a counterexample to the Zassenhaus 
conjecture. One can easily compute its outer derivation algebra over $\F_3$ and hence over $\F$. 
This is well known in the literature, see for example \cite{BEC}, Theorem $1$, but hasn't been associated explicitly 
to the Zassenhaus conjecture.

\begin{prop}\label{3.4}
Let $\F$ be a field of characteristic $p=3$. Then the derivation algebra of $\Lg=\mathfrak{psl}_3(\F)$ is given by
$\Der(\Lg)\cong \Lg_2$, and the quotient by $\Inn(\Lg)\cong \Lg$ is given by $\Out(\Lg)\cong \Lg$. In particular
the outer derivation algebra of $\Lg$ is simple and non-solvable.
\end{prop}

For the series of simple Lie algebras $\mathfrak{psl}_{3m}(\F)$ over a field $\F$ of characteristic $3$ 
the case $m=1$ is the only exception regarding Proposition $\ref{3.3}$. We have the following result, which follows from
Corollary $3.4$ on page $73$ in \cite{ELK}.

\begin{prop}\label{3.5}
Let $\F$ be an algebraically closed field of characteristic $3$. Then we have  
$\Der (\mathfrak{psl}_{3m}(\F))\cong \mathfrak{pgl}_{3m}(\F)$ for all $m\ge 2$. Hence $\Out(\mathfrak{psl}_{3m}(\F))\cong \F$ 
is abelian for all $m\ge 2$. 
\end{prop}

For $\Le_6$ we have the following result.

\begin{prop}
Let $\F$ be a field of characteristic $p=3$ and $\Lg=\Le_6/Z$ be the simple Lie algebra of dimension $77$. 
Then $\Out(\Lg)$ is abelian.
\end{prop}

\begin{proof}
In fact, $\Lg$ has a derivation algebra of dimension $78$, so that the outer derivation algebra is $1$-dimensional. 
This follows from a computation with GAP. 
\end{proof}

For $p=2$ the only simple Chevalley algebras are $A_{2n}, G_2,E_6,E_8$. The algebras of type 
$A_{2n+1},B_n,C_n,D_n,E_7$ have a non-zero center, and $F_4$ has a $26$-dimensional ideal $J$ generated by the 
short roots. We have $F_4/J\cong J$ and this is a simple Lie algebra in characteristic $2$, see \cite{BJG}, or
\cite{PUR}, section $5.1$. Furthermore $F_4$ has up to isomorphism only one such simple ideal of dimension $26$.
This yields another counterexample to the Zassenhaus conjecture, see \cite{BJG}. 

\begin{prop}\label{3.7}
Let $\F$ be a field of characteristic $p=2$ and $J$ be the simple ideal in $F_4$
generated by the short roots. Then $\Der(J)\cong F_4$ and $\Out(J)\cong J$ is simple and non-solvable.
\end{prop}

According to Seligman \cite{SEL}, this has been first shown in \cite{SCT}, for an algebraically closed 
field of characteristic $2$.

\subsection{Cartan  type}
The second class of simple modular Lie algebras consists of the Lie algebras of {\em Cartan type}
\[
W(m;\underline{n}),\; S(m;\underline{n};\Phi)^{(1)},\; H(2r;\underline{n};\Phi)^{(2)},\; K(2r+1;\underline{n})^{(1)}.
\]
The characteristic zero versions of these Lie algebras are infinite dimensional Lie algebras occurring in E. Cartan's work
of $1909$ on pseudogroups in differential geometry. For the precise definition see \cite{STR1}, Chapter $4$.
These Lie algebras are filtered deformations of a special class, namely of the {\em graded} simple algebras of Cartan type, 
denoted by $W(m;\underline{n}),\; S(m;\underline{n})^{(1)},\; H(2r;\underline{n})^{(2)},\; K(2r+1;\underline{n})^{(1)}$.
Let $\underline{n}=(n_1,\ldots ,n_m)$ and $s=n_1+\cdots +n_m$. The {\em Witt algebra} $W(m;\underline{n})$ has 
dimension $mp^s$ and is simple unless $p=2$ and $m=1$. The {\em special algebra} $S(m;\underline{n})^{(1)}$ 
has dimension $(m-1)(p^s-1)$ and is simple for all $p>0$ and all $m\ge 1$. 
The {\em Hamiltonian algebra} $H(2r;\underline{n})^{(2)}$ is simple of
dimension $p^s-2$, and finally the {\em contact algebra} $K(2r+1;\underline{n})^{(1)}$ is simple with
\[
\dim (K(2r+1;\underline{n})^{(1)})=\begin{cases} p^s,  \hspace*{0.7cm} \text{ if } 2r+1\not\equiv -3 \bmod p,\\
p^s-1, \text{ if } 2r+1\equiv -3 \bmod p.\end{cases}
\]
For $p=3$ or $p=2$ some of these algebras coincide with algebras of classical type, e.g., we have
$H(2;\underline{1})^{(2)}\cong \mathfrak{psl}(3)$ with $\underline{1}=(1,1)$ for $p=3$. This does not happen for $p>3$. \\[0.2cm]
We have the following result.

\begin{prop}\label{3.8}
Let $\Lg$ be a simple Lie algebra of Cartan type over an algebraically closed field of characteristic $p>3$.
Then $\Out(\Lg)$ is solvable. More precisely, $\Out(\Lg)$ is solvable of derived length $d\le 1$ for type $W$ and type
$K$, of derived length $d\le 2$ for type $S$ and of derived length $d\le 3$ for type $H$. 
\end{prop}

\begin{proof}
The proof is due to Leinen and Puglisi in \cite{LEP}, Proposition $5$. They showed that $\Out(\Lg)$ is solvable of 
derived length $d\le 4$ over an algebraically closed field of characteristic $p>7$, using results of \cite{BKK}.
Here the assumption $p>7$ was made because the classification at this time was only known for $p>7$.
However, all results used are actually stated in \cite{BKK} for $p>3$. Hence we immediately obtain the result
of Proposition $5$ for all $p>3$. It remains to refine the proof, so that we obtain derived length $d\le 3$. We will
use the notation of \cite{BKK}. Let $\ad_{\Lg}(\Lh)=\{\ad(x)_{\mid \Lg}\mid x\in \Lh\}$ for a Lie algebra $\Lh$ and an ideal $\Lg$ 
in $\Lh$. \\[0.2cm]
First let $\Lg=W(m;\mathbf{n})$ be the Witt algebra. Consider its $p$-envelope $\ov{\Lg}$ in the derivation
algebra $\Der(A(m;\mathbf{n}))$. Then we have $\Der(\Lg)=\ad_{\Lg}(\ov{\Lg})$, see \cite{BKK} p. $903$.
Hence it follows that
\[
\Der(\Lg)^{(1)}=\ad_{\Lg}(\ov{\Lg})^{(1)}=\ad_{\Lg}(\ov{\Lg}^{(1)})=\ad_{\Lg}(\Lg^{(1)})=\Inn(\Lg).
\]
Here we have used that $\overline{\Lg}^{(1)}=\Lg^{(1)}$. Hence $\Out(\Lg)$ is abelian. The same proof
works for the contact Lie algebra $\Lg=K(m;\mathbf{n})$. \\[0.2cm]
Let $S=S(m;\mathbf{n},\om)$ be the special Lie algebra. There are inclusions 
$S\subseteq \widetilde{S}\subseteq C\widetilde{S}$ of algebras in $\Der(A(m;\mathbf{n}))$, see
\cite{BKK}, section $1$ for the definitions. Let $P$ denote the $p$-envelope of $C\widetilde{S}$ in $\Der(A(m;\mathbf{n}))$.
Then we have $\Der(S)=\ad_S(P)$ by \cite{BKK}, p. $903$ and hence
\[
\Der(S)^{(1)}=\ad_S(P)^{(1)}=\ad_S(C\widetilde{S}^{(1)}).
\]
Since $\widetilde{S}$ is an ideal of $C\widetilde{S}$ of codimension at most one, the quotient 
$C\widetilde{S}/\widetilde{S}$ is abelian and we have $(C\widetilde{S})^{(1)}\subseteq \widetilde{S}$.
It follows that $\Der(S)^{(1)}\subseteq \ad_S(\widetilde{S})$ and $\Der(S)^{(2)}\subseteq \ad_S(\widetilde{S}^{(1)})=\Inn(S)$.
Hence $\Out(S)$ is solvable of derived length at most two. \\[0.2cm]
Let $H=H(m;\mathbf{n},\om)$ be the Hamiltonian Lie algebra. Similarly to type $S$, it follows from \cite{BKK} that
$\Der(H)^{(3)}\subseteq \ad_S(\widetilde{H}^{(2)})=\Inn(H)$. Hence $\Out(H)$ is solvable of derived length at most three. 
\end{proof}

It is easy to see that the bounds are sharp. In particular, there are simple Lie algebras of Cartan type where
the outer derivation algebra has exactly derived length $d=3$ and hence is not metabelian.

\begin{ex}
Let $\F$ be a field of characteristic $p>3$ and $\Lg=H(2;\underline{1})^{(2)}$ be the Hamiltonian algebra of dimension 
$p^2-2$. Then $\Out(\Lg)$ is solvable of derived length $d=3$. 
\end{ex}

This has been observed by J. Feldvoss. In fact $\Out(\Lg)^{(1)}$ is isomorphic to the $3$-dimensional Heisenberg Lie algebra.

\subsection{Melikian  type}
In characteristic $p=5$ another family of simple Lie algebras is known, different from the classical or Cartan type.
The algebras $\CM (n_1,n_2)$ are called {\em Melikian algebras} with parameters $(n_1,n_2)\in \N^2$. They are simple
of dimension $5^{n_1+n_2+1}$. 

\begin{prop}\label{3.10}
Let $\CM$ be a Melikian algebra over an algebraically closed field of characteristic $5$. Then $\Out(\CM)$ is abelian.
\end{prop}

\begin{proof}
The proof works exactly the same way as the proof of Proposition $\ref{3.8}$ for type $W$ or $K$ because we have 
$\Der(\CM)=\ad_{\CM}(\ov{\CM})$, see Theorem $7.1.4$ and Theorem $7.2.2$ in \cite{STR1}.
\end{proof}

\subsection{Summary}
Let us summarize the results on the Zassenhaus conjecture for algebraically closed fields of characteristic $p>3$ and
add some results for characteristic $2$ and $3$. By Proposition $\ref{3.3}$, $\ref{3.8}$ and $\ref{3.10}$ the Zassenhaus
conjecture has a positive answer in the following case.

\begin{thm}
Let $\Lg$ be a simple modular Lie algebra over an algebraically closed field of characteristic $p>3$. Then
$\Out(\Lg)$ is solvable of derived length at most three.
\end{thm}

If $\Lg$ is a central simple Lie algebra over an arbitrary field $\F$ of characteristic $p>3$, then $\Lg \otimes_{\F}\ov{\F}$
is simple over $\ov{\F}$. Hence we obtain the following corollary.

\begin{cor}\label{3.12}
Let $\Lg$ be a central simple Lie algebra over a field of characteristic $p>3$. Then
$\Out(\Lg)$ is solvable of derived length at most three.
\end{cor}

\begin{rem}
As we have seen in Propositions $\ref{3.7}$ and $\ref{3.4}$ there are counterexamples to the Zassenhaus conjecture 
for $p=2$ and $p=3$ arising from the simple constituents of the classical Lie algebras $\Lf_4$ and $\Lg_2$.
On the other hand, we haven't found other counterexamples so far among the known new simple Lie algebras 
over an algebraically closed field in characteristic $3$ and $2$. In \cite{CNZ}, published in $1993$, the outer 
derivation algebras for several families of new simple Lie algebras in characteristic $3$ are computed. It is noted that it 
is always solvable. In particular, all derivations are inner for the algebras $L(\epsilon)$ of Kostrikin, and 
$\Out(\Lg)$ is abelian for the algebras of type $R$ and $T$ of Kuznetsov and the algebras $Y(\underline{n})$ of 
Skryabin. \\[0.2cm]
For characteristic $p=2$ there are apparently too many simple Lie algebras to obtain conclusive results
and we only have verified the Zassenhaus conjecture for lists of known simple Lie algebras in low dimensions over $\F_2$. 
In particular, we have studied the lists of
M. Vaughan-Lee \cite{VAU} and B. Eick \cite{EIK}. The result is that all algebras have a solvable outer derivation 
algebra of derived length at most three. In the list of \cite{EIK}, up to dimension $n\le 20$, there are only 
two different algebras, namely $\Lg=S(2,2)$ and ${\rm Kap}_1(4)$ of dimension $14$, where $\Out (\Lg)$ has
derived length $d=3$.
\end{rem}

\section{Taut Lie algebras}

For our results concerning CPA-structures on simple modular Lie algebras we will need a short term for the following
property of a Lie algebra.

\begin{defi}
Let $\Lg$ be a Lie algebra. We say that $\Lg$ is {\em taut} if every Lie algebra homomorphism
$\phi\colon \Lg\ra \Der(\Lg)$ maps $\Lg$ into $\Inn(\Lg)$.
\end{defi}

Clearly every Lie algebra whose derivations are all inner is taut. However, the derivations of a simple
modular Lie algebra need not be inner in general. But we have the following result.

\begin{prop}\label{4.2}
Let $\Lg$ be a simple Lie algebra with $\dim (\Out (\Lg))<\dim (\Lg)$. Then $\Lg$ is taut.
\end{prop}

\begin{proof}
Let $\phi\colon \Lg\ra \Der(\Lg)$ be an arbitrary Lie algebra homomorphism and let $\pi\colon \Der(\Lg)\ra \Out(\Lg)$ 
be the natural projection. Since $\Lg$ is simple, $\phi(\Lg)$ is a subalgebra of $\Out(\Lg)$, which is either trivial 
or isomorphic to $\Lg$. But the latter cannot happen because of
$\dim (\Lg)>\dim (\Out (\Lg))$. Hence $\phi(\Lg)\subseteq \ker(\pi)=\Inn(\Lg)$, so that $\phi$ is taut. 
\end{proof}

Taut Lie algebras are related to the Zassenhaus conjecture as follows.

\begin{prop}\label{4.3}
Let $\Lg$ be a simple modular Lie algebra satisfying the Zassenhaus conjecture. Then $\Lg$ is taut.
\end{prop}

\begin{proof}
Recall that $\Lh^{(\infty)}=\cap_{i\in \N}\Lh^{(i)}$ for a Lie algebra $\Lh$. Let  $\phi\colon \Lg\ra \Der(\Lg)$ be a Lie 
algebra homomorphism. Then we have
\[
\phi(\Lg)=\phi(\Lg^{(\infty)})=\phi(\Lg)^{(\infty)}\subseteq \Der(\Lg)^{(\infty)}\subseteq \Inn(\Lg)
\]
by the assumptions.
\end{proof}

\begin{cor}\label{4.4}
Let $\Lg$ be a central simple Lie algebra over a field of characteristic $p>3$. Then $\Lg$ is taut.
\end{cor}

\begin{proof}
By Corollary $\ref{3.12}$ we know that $\Out(\Lg)$ is solvable, so that the claim follows from Proposition $\ref{4.3}$.
\end{proof}

Is there a simple modular Lie algebra which is not taut? If so, it must be a counterexample to the Zassenhaus
conjecture. In fact, for $p=2$ and $p=3$ we have such a counterexample, see Proposition $\ref{3.7}$ and Proposition
$\ref{3.4}$. These two simple Lie algebras $\Lg$ satisfy $\Lg\cong \Out(\Lg)$. We have the following lemma.

\begin{lem}\label{4.5}
Let $\Lg$ be a simple Lie algebra satisfying $\Lg\cong \Out(\Lg)$. Then $\Lg$ is taut if and only if the short exact
sequence
\[
0\ra \Inn(\Lg) \ra \Der(\Lg)\ra \Out(\Lg)\ra 0
\]
does not split.
\end{lem}

\begin{proof}
We will show that $\Lg$ is not taut if and only if the above short exact sequence splits. Hence suppose first that
$\Lg$ is not taut, i.e., there is a homomorphism $\phi\colon \Lg\ra \Der(\Lg)$ with $\phi(\Lg)$ not contained
in $\Inn(\Lg)$. Denote by $\pi\colon \Der(\Lg)\ra \Out(\Lg)$ the canonical projection. Then $\phi(\Lg)$ is
nonzero. Since $\Lg$ is simple we have $\phi(\Lg)\cong \Lg$. Then the ideal $\Inn(\Lg)\cap \phi(\Lg)$ is trivial.
So the subalgebra $\Inn(\Lg)\rtimes \phi(\Lg)$ of $\Der(\Lg)$ induces a split extension. Since their dimensions
coincide, we conclude that $\Der(\Lg)$ is a split extension of $\Inn(\Lg)$ by $\phi(\Lg)$. \\
Conversely assume that the above sequence splits. Then there exists an embedding $\phi\colon \Out(\Lg)\ra \Der(\Lg)$
such that $\phi(\Out(\Lg))$ is not contained in $\Inn(\Lg)$. Let $\theta\colon \Lg\ra \Out(\Lg)$ be an isomorphism.
Then $(\theta\circ \phi)(\Lg)$ is not contained in $\Inn(\Lg)$, so that $\Lg$ is not taut.
\end{proof}

For the counterexample given in Proposition $\ref{3.4}$ we have the following result over $\F_3$.

\begin{prop}\label{4.6}
The Lie algebra $\mathfrak{psl}_3(\F_3)$ is taut.
\end{prop}

\begin{proof}
Let  $\Ls=\mathfrak{psl}_3(\F_3)$. We have $\Inn(\Ls)\cong \Out(\Ls)\cong \Ls$ and $\Der(\Ls)\cong \Lg_2$, 
see Proposition $\ref{3.4}$. We will show by a direct computation that the short exact sequence
\[
0\ra \Ls \ra \Lg_2 \xrightarrow{\pi} \Ls \ra 0
\]
does not split. Then the claim follows from Lemma $\ref{3.5}$. Let  
\begin{align*}
e_1 & =E_{12},\, e_2=E_{13},\, e_3=E_{21},\, e_4=E_{23},\, e_5=E_{31},\, e_6=E_{32},\\
e_7 & =E_{11}-E_{22},\, e_9=E_{22}-E_{33}
\end{align*}
be a basis for $\mathfrak{sl}_3(\F_3)$. Then the center is generated by $I_3=e_7-e_8$.
This yields a basis for $\Ls$, again denoted by $(e_1,\ldots ,e_7)$ with Lie brackets
\begin{align*}
[e_1,e_3] & = e_7,  & [e_2,e_5] & =2e_7, & [e_4,e_5]& =e_3,\\
[e_1,e_4] & =e_2,   & [e_2,e_6] & =e_1 , & [e_4,e_6]& =e_7, \\
[e_1,e_5] & = -e_6, & [e_2,e_7] & =-e_2, & [e_4,e_7]& =e_4, \\
[e_1,e_7] & =-2e_1, & [e_3,e_6] & =-e_5, & [e_5,e_7]& =e_5, \\
[e_2,e_3] & =-e_4,  & [e_3,e_7] & =2e_3, & [e_6,e_7]& =-e_6.
\end{align*}
Suppose that the sequence splits. Then there is a Lie algebra homomorphism $\phi \colon \Ls\ra \Lg_2$ such that 
$\pi\circ \phi=\id_{|\Ls}$. So $\phi$ has a $14\times 7$-matrix of block form 
\[
\phi=
\left(
\begin{array}{c}
A \\
\hline
I
\end{array}
\right)
\]
with the identity matrix $I=I_7$ and $A=(a_{ij})\in M_7(\F_3)$. The condition that $\phi$ is a homomorphism
is equivalent to a system of polynomial equations in $49$ the variables $a_{ij}$ over $\F_3$. Specifying a particular value for 
a column vector $v\in \F_3^7$ of $A$, it can be easily solved using Gr\"obner bases over $\F_3$. It is no problem to check 
all possible $3^7$ cases to see that there is no solution.
\end{proof}

\begin{rem}
We also believe that the $26$-dimensional simple Lie algebra $\Ls$ over $\F_2$ of Proposition $3.7$ is taut, because again 
the short exact sequence 
\[
0\ra \Ls \ra \Lf_4 \xrightarrow{\pi} \Ls \ra 0
\]
probably does not split. However, we were not able to solve the $2^{26}$ systems of polynomial equations in $676$ variables over
$\F_2$. 
\end{rem}

Furthermore we have checked several other families of simple Lie algebras in characteristic $2$ and $3$ 
and haven't found an example which is not taut. So we have the following conjecture.

\begin{con}
Every simple modular Lie algebra is taut.
\end{con}

\section{CPA-structures on simple modular Lie algebras}

We recall the following definition from \cite{BU52}.

\begin{defi}
A CPA-structure on a Lie algebra $\Lg$ is called {\em inner}, if there is a Lie algebra endomorphism
$\phi\in \End(\Lg)$ such that the CPA-product is given by $x\cdot y=[\phi(x),y]$.
\end{defi}

If $\Lg$ has trivial center, then any linear map $\phi\in \End(V)$ defining a CPA-structure this way is already
a Lie algebra homomorphism, see Lemma $2.8$ in \cite{BU52}. It is easy to see that we can generalize many structure
results of \cite{BU52} on inner CPA-structures to fields of arbitrary characteristic. In particular Lemma $2.11$
of \cite{BU52} generalizes as follows.

\begin{lem}\label{5.2}
Let $\Lg$ be a Lie algebra over an algebraically closed field $\F$ of arbitrary characteristic and suppose that 
$x\cdot y=[\phi(x),y]$ defines an inner CPA-structure on $\Lg$. Let
$\Lg=\bigoplus_{\al}\Lg_{\al}$ be the generalized eigenspace decomposition of $\Lg$ with respect to $\phi$. Then we have
\[
[\Lg_{\al},\Lg_{\be}]  \subseteq \Lg_{-\al^2} \cap \Lg_{\al\be}\cap \Lg_{-\be^2}
\]
for all $\al,\be \in \F$. In particular, $[\Lg_{\al},\Lg_{\be}] \neq 0$ implies $\al+\be=0$.
\end{lem}

For elements $x_1,\ldots ,x_n$ in a Lie algebra $\Lg$ define 
\begin{align*}
[x_1,x_2,x_3]& = [x_1,[x_2,x_3]],\\
[x_1,x_2,x_3,x_4]& =[x_1,[x_2,[x_3,x_4]]]
\end{align*}
and so forth. By induction we obtain the following commutator formula.

\begin{lem}\label{5.3}
Let $\Lg$ be a Lie algebra and suppose that $x\cdot y=[\phi(x),y]$ defines an inner CPA-structure on $\Lg$. Then
for all $n\ge 1$, for all elements $x_1,\ldots ,x_n$ and $y$ in $\Lg$ we have
\[
[x_1,\ldots ,x_n]\cdot y =\pm [[x_1,\ldots ,x_n],\phi^{f(n)}(y)],
\]
where the map $f\colon \N\ra \N$ is defined recursively by $f(1)=1$ and $f(n+1)=2f(n)+1$.
\end{lem}

We have the following structure result on inner CPA-structures.

\begin{prop}\label{5.4}
Let $\Lg$ be a Lie algebra over an algebraically closed field $\F$ of arbitrary characteristic $p\ge 0$ and suppose 
that $x\cdot y=[\phi(x),y]$ defines an inner CPA-structure on $\Lg$. Let $\Lg=\bigoplus_{\al\in \F}\Lg_{\al}$ be the generalized 
eigenspace decomposition of $\Lg$ with respect to $\phi$. Then $\Lg$ decomposes into a direct sum 
of ideals 
\[
\Lg=\Lg_0\oplus \bigoplus_{\al\neq 0}\Lg_{\al}=\Ln\oplus \Lh 
\]
with the following properties. 
\vspace{0.1cm}
\begin{itemize}
\item[(1)] $(v,v')\cdot (w,w')=(v\cdot w, v'\cdot w')$ for all $v,w\in \Ln$ and $v',w'\in \Lh$. 
\vspace{0.1cm}
\item[(2)] $\Ln^{\infty}\cdot \Ln=0$. 
\vspace{0.1cm}
\item[(3)] If $p\neq 2$ then $[[\Lh,\Lh],[\Lh,\Lh]]=0$. 
\vspace{0.1cm}
\item[(4)] If $p=2$ then $v\cdot w=[v,w]$ for all $v\in [\Lh,[\Lh,\Lh]]$ and $w\in \Lh$.
\end{itemize}
\end{prop}

\begin{proof}
We follow the proof of Theorem $2.1.4$ in \cite{BU52}. For $\al\in \F$ we have $[\Lg_0,\Lg_{\al}]\subseteq \Lg_0$, 
so that $\Ln$ is an ideal of $\Lg$. For all $\al,\be \in \F^{\times}$ we have $[\Lg_{\al},\Lg_{\be}]\subseteq \Lg_{\al\be}$
and $[\Lg_{0},\Lg_{\be}]\subseteq \Lg_{-\be^2}$, so that $\Lh$ is an ideal of $\Lg$. By Lemma $2.1$ of \cite{BU52}, $\Ln$
and $\Lh$ are also ideals for the CPA-structure. By construction we have $\phi(\Ln)\subseteq \Ln$ and 
 $\phi(\Lh)\subseteq \Lh$, so that we may look at the restrictions $\phi_{\Ln}$ and $\phi_{\Lh}$ of $\phi$ to
$\Ln$ respectively $\Lh$. For $(x,y)$ and $(x',y')$ in $\Lg=\Ln\oplus \Lh$ the CPA-structure is given by the formula
$(x,y)\cdot (x',y')=([\phi_{\Ln}(x),x'],[\phi_{\Lh}(y),y'])$. This shows $(1)$. By definition the endomorphism $\phi_{\Ln}$
is nilpotent so that the commutator formula $\ref{5.3}$ implies $(2)$. As in the proof of Theorem $2.1.4$ in \cite{BU52}
we conclude that we either have $\Lh^{(2)}=[[\Lh,\Lh],[\Lh,\Lh]]=0$ or 
\[
\Lh^{(2)} \subseteq (\Lg_i\cup \Lg_{-i})\cap (\Lg_1\cup \Lg_{-1}),
\]
where $\al=i$ is a root of the polynomial $X^2+1$ in $\F$. In the second case we have 
\[
\{i,-i\}\cap \{1,-1\}=\emptyset
\]
for $p\neq 2$, so that again $\Lh^{(2)}=0$. This shows $(3)$. Assume now that $p=2$. Then we have $i=-i=1=-1$ so that
$\Lh^{(2)} \subseteq \Lg_1$. Hence by construction the restriction of $\phi$ to $\Lh^{(2)}$ is unipotent. Let
$x_1,x_2,y \in \Lh$. By Lemma $\ref{5.3}$ we have
\begin{align*}
[[x_1,x_2],\phi_{\Lh}(y)] & = [\phi_{\Lh}(y),[x_1,x_2]] \\
                         & = y\cdot [x_1,x_2] \\
                         & = [x_1,x_2]\cdot y \\
                         & = [[x_1,x_2],\phi_{\Lh}^3(y)].
\end{align*}
Since $\phi_{\Lh}$ is an automorphism of $\Lh$ this implies 
\[
[[\Lg,\Lg],(\phi_{\Lh}^2-\id_{\Lh})(\Lg)]\subseteq 0.
\]
Now let $x\in [\Lh,\Lh]$ and $y,z \in \Lh$. Then 
\[
[\phi_{\Lh}(z),[x,y]]=[z,[\phi_{\Lh}(x),\phi_{\Lh}(y)]]=[z,[x,\phi_{\Lh}^2(y)]]=[z,[x,y]].
\]
So we have indeed $v\cdot w=[v,w]$ for all $v\in [\Lh,[\Lh,\Lh]]$ and $w\in \Lh$.
\end{proof}

Now we are able to prove results for CPA-structures on simple modular Lie algebras.

\begin{prop}\label{5.5}
Let $\Lg$ be a Lie algebra over a field $\F$ of characteristic $p\neq 2$. If $\Lg$ is taut, then every CPA-structure on
$\Lg$ satisfies $\Lg^{(\infty)}\cdot \Lg = 0$.
\end{prop}

\begin{proof}
The left multiplication $L\colon \Lg\ra \Der(\Lg)$, $x\mapsto L(x)$ is a Lie algebra homomorphism. Since 
$\Lg$ is taut, it follows that $L(x)\in \Inn(\Lg)$ for all $x\in \Lg$. Hence there is a linear map
$\phi\colon \Lg\ra \Lg$ such that $L(x)=\ad(\phi(x))$, i.e., $x\cdot y=[\phi(x),y]$ for all $x,y\in \Lg$.
By extending the base field $\F$ to an algebraic closure $\ov{\F}$ of it we can naturally extend $\Lg$ to 
$\ov{\Lg}=\F\otimes_{\F}\Lg$ and the CPA-structure $L\colon \Lg\ra \Inn(\Lg)$ to a CPA-structure 
$\ov{L}\colon \ov{\Lg}\ra \Inn(\ov{\Lg})$ by  $1\otimes v\mapsto 1\otimes L(v)$. Let $\Lk:=\ov{\Lg}/Z_{\infty}(\ov{\Lg})$,
where $Z_{\infty}(\ov{\Lg})$ is the last term of the upper-central series of $\ov{\Lg}$.
Then $Z(\Lk)=0$, so that the induced CPA-structure $\widetilde{L}\colon \Lk\ra \Lk$ is inner by Lemma $2.8$ in \cite{BU52}. 
This says that $\phi$ is a Lie algebra endomorphism.
By Proposition $\ref{5.4}$, $(2),(3)$ it follows that $\Lk^{(\infty)}\cdot \Lk \subseteq (\Lk^{(2)})^{\infty}\cdot \Lk=0$, so 
that $\ov{\Lg}^{(\infty)}\cdot \ov{\Lg}\subseteq Z_{\infty}(\ov{\Lg})$. Since
$[\ov{\Lg}^{(\infty)},\ov{\Lg}^{(\infty)}]=\ov{\Lg}^{(\infty)}$ we also have
\[
\ov{\Lg}^{(\infty)}\cdot \ov{\Lg}\subseteq \ov{\Lg}^{(\infty)}\cdot (\ov{\Lg}^{(\infty)}\cdot \ov{\Lg})\subseteq 
\ov{\Lg}^{(\infty)}\cdot Z_{\infty}(\ov{\Lg}) =0.
\]
But this implies that $\Lg^{(\infty)}\cdot \Lg = 0$.
\end{proof}

\begin{cor}\label{5.6}
Let $\Lg$ be a perfect Lie algebra over a field of characteristic $p\neq 2$. Suppose that $\Lg$ is taut. 
Then every CPA-structure on $\Lg$ is trivial. In particular, every perfect Lie algebra over such a field
having a solvable outer derivation algebra only admits the trivial CPA-structure. 
\end{cor}

\begin{proof}
Since $\Lg$ is perfect we have $\Lg^{(\infty)}=\Lg$. If $\Lg$ has a solvable outer derivation algebra, then it is taut by
Proposition $\ref{4.3}$. Hence we have $\Lg\cdot \Lg=0$ by Proposition $\ref{5.5}$.
\end{proof}

We come to a main result.

\begin{thm}\label{5.7}
Let $\Lg$ be a central simple Lie algebra of characteristic $p>3$. Then every CPA-structure on $\Lg$ is trivial.
\end{thm}

\begin{proof}
Since $\Lg$ is taut by Corollary $\ref{4.4}$ and also perfect the claim follows from Corollary $\ref{5.6}$.
\end{proof}

What can we say about CPA-structures of simple modular Lie algebras in characteristic $2$ or $3$? Let us first investigate
the case $p=3$. All known simple modular Lie algebras in characteristic $3$ only admit the trivial CPA-structure, it seems.
We just give one example here.

\begin{ex}
Every CPA-structure on $\mathfrak{psl}_3(\F_3)$ is trivial.
\end{ex}

This follows either by an easy direct computation, or by Corollary $\ref{5.6}$. Indeed, $\mathfrak{psl}_3(\F_3)$
is perfect and it is taut by Proposition $\ref{4.6}$. Recall that a simple Lie algebra with solvable outer derivation algebra
is taut, but that the converse need not be true, as we see here. \\[0.2cm]
In characteristic $2$ the Lie bracket is commutative, so that $x\cdot y=[x,y]$ defines a non-trivial CPA-structure on any
non-trivial Lie algebra. We call this structure the {\em adjoint} CPA-structure.
We can show the following result.

\begin{prop}\label{5.9}
Let $\Lg$ be a simple modular Lie algebra in characteristic $2$. If $\Lg$ is taut, then every CPA-structure on $\Lg$ is either
trivial or adjoint.
\end{prop}

\begin{proof}
Let $L\colon \Lg\ra \Der(\Lg)$ be the homomorphism induced by left multiplication. Since $\Lg$ is taut we have
$L(\Lg)\subseteq \Inn(\Lg)$. Since $Z(\Lg)=0$ the CPA-structure on $\Lg$ is inner, i.e., there exists
a unique homomorphism $\phi\colon \Lg\ra \Lg$ such that $x\cdot y=L(x)(y)=[\phi(x),y]$. Since $\Lg$ is simple,
$\phi(\Lg)$ is either trivial or isomorphic to $\Lg$. In the first case we have $\Lg\cdot \Lg=0$. So we may
assume that $\phi$ is an isomorphism. Let $\Lk$ be the perfect Lie algebra obtained from $\Lg$ by extending
the scalars over an algebraic closure of the ground field. Then $\phi$ naturally extends to an isomorphism
$\psi\colon \Lk\ra \Lk$ and also defines an inner CPA-structure on $\Lk$.  Now we can apply  part $(4)$ of 
Proposition $\ref{5.4}$ and decompose $\Lk=\Ln\oplus \Lh$ as a direct sum of Lie algebras. Since $\psi$ is
an isomorphism, we have $\Ln=0$ and $\Lk=\Lh$. Since $\Lk$ is perfect we have $\psi=\id_{\mid \Lk}$.
By restriction we obtain $\phi=\psi_{\mid \Lg}=\id_{\Lg}$, so that $x\cdot y=[x,y]$ for all $x,y\in \Lg$.
\end{proof}

\begin{ex}
Let $\Lg$ be a simple modular Lie algebra over $\F_2$ from the list of B. Eick \cite{EIK}. Then every CPA-structure on
$\Lg$ is trivial or adjoint.
\end{ex}

Indeed, every such algebra is taut. \\[0.2cm]
The following example shows that the CPA-structures may also depend on the base field.

\begin{ex}
Let $\Lg=W(1,\underline{2})^{(1)}\otimes_{\F_2}\F_4$ be the simple modular Lie algebra over $\F_2$ from Example $\ref{2.3}$.
Then $\Lg$ admits only the trivial and adjoint CPA-structure over $\F_2$, but admits other CPA-structures over an
algebraic closure of $\F_2$.
\end{ex}

One can explain how the other CPA-structures arise. In fact $\Lg$ is not central simple. It decomposes over
$\ov{\F}_2$ as the sum $\Lg=\La\oplus \Lb$ of two simple ideals. There is a basis $(x_1,\ldots ,x_6)$ of $\Lg$ such that 
the Lie brackets are given by
\begin{align*}
[x_1,x_2] & = x_3,      & [x_2,x_3] & =x_1,     & [x_3,x_6]& =x_2,\\
[x_1,x_3] & = x_4,      & [x_2,x_5] & =x_1+x_6, & [x_4,x_5]& =x_6, \\
[x_1,x_4] & = x_5,      & [x_2,x_6] & =x_3+x_5, & [x_4,x_6]& =x_3, \\
[x_1,x_5] & = x_2+x_4,  & [x_3,x_4] & =x_1+x_6, & [x_5,x_6]& =x_4.
\end{align*}
All CPA-structures on $\Lg$ are explicitly given as follows, where $\al\in \ov{\F}_2$ is a root of 
\[
X(X^3-1)=X(X-1)(X^2+X+1),
\]
\begin{align*}
x_1\cdot x_2 & = \al x_3 +\al(\al+1)x_5,  & x_2\cdot x_3 & =\al^2 x_1+ \al(\al+1)x_6, & x_3\cdot x_6 & =\al x_2+ \al(\al-1)x_4,\\
x_1\cdot x_3 & = \al(\al+1)x_2+\al^2x_4,  & x_2\cdot x_5 & =\al x_1+\al^2 x_6,        & x_4\cdot x_5 & =\al(\al+1)x_1+\al x_6, \\
x_1\cdot x_4 & = \al(\al+1)x_3+ \al^2x_5, & x_2\cdot x_6 & =\al^2 x_3+\al x_5, & x_4\cdot x_6 & =\al x_3+\al(\al+1)x_5, \\
x_1\cdot x_5 & = \al^2x_2+\al x_4,        & x_3\cdot x_4 & =\al x_1+\al^2 x_6, & x_5\cdot x_6 & =\al(\al+1)x_2+\al^2 x_4.
\end{align*}
For $\al=0$ we obtain the trivial CPA-structure and for $\al=1$ the adjoint one. The remaining ones are given by $\al$
with $\al^2+\al+1=0$ over $\overline{\F}_2$. Altogether we just obtain the four CPA-structures from the decomposition
$\Lg=\La\oplus \Lb$, namely $\ad_{\La}\oplus \ad_{\Lb}$, $\ad_{\La}\oplus 0_{\Lb}$, $0_{\La}\oplus \ad_{\Lb}$ and
$0_{\La}\oplus 0_{\Lb}$.

\begin{rem}
It is also possible to generalize results on post-Lie algebra structures on pairs of Lie algebras from complex
numbers to algebraically closed fields of characteristic $p>3$. For example, one can generalize Proposition $4.6$
of \cite{BU41} as follows: every post-Lie algebra structure on a pair $(\Lg,\Ln)$ with $\Lg$ and $\Ln$ simple, over
an algebraically closed field of characteristic $p>3$ satisfies either $x\cdot y=0$ or $x\cdot y=[x,y]=-\{x,y\}$,
where the first Lie bracket is in $\Lg$ and the second one in $\Ln$. The proof over the complex numbers uses
algebraic arguments which are still valid, except for the fact that $\Out(\Ln)$ is no longer trivial. However the proof
still works by using that $\Out(\Ln)$ is solvable in this case.
\end{rem}

\section*{Acknowledgments}
Dietrich Burde is supported by the Austrian Science Foun\-da\-tion FWF, grant I3248. 
Wolfgang A. Moens acknowledges support by the Austrian Science Foun\-da\-tion FWF, 
grant P30842. We thank Bettina Eick, J\"org Feldvoss and Willem de Graaf for helpful discussions.

\end{document}